\documentclass[aip,
 amsmath,amssymb,
 reprint,%
]{revtex4-1}

\usepackage{graphicx}% Include figure files
\usepackage{dcolumn}% Align table columns on decimal point
\usepackage{bm}% bold math
\usepackage[utf8]{inputenc}
\usepackage[T1]{fontenc}
\usepackage{mathptmx}
\usepackage{etoolbox}

\makeatletter
\def\@email#1#2{%
 \endgroup
 \patchcmd{\titleblock@produce}
  {\frontmatter@RRAPformat}
  {\frontmatter@RRAPformat{\produce@RRAP{*#1\href{mailto:#2}{#2}}}\frontmatter@RRAPformat}
  {}{}
}%
\makeatother
\begin{document}

\preprint{AIP/123-QED}

\title{Mixed Mode Bursting Oscillations Induced by Birhythmicity and Noise}
% Force line breaks with \\
\author{Na Yu}
\email{nayu@ryerson.ca}
\affiliation{Department of Mathematics, Ryerson University}
\author{Xuan Xia}%
 \affiliation{Department of Mathematics, Ryerson University}
 \author{Juan Liyau}
 \affiliation{Department of Mathematics, Ryerson University}

            
\begin{abstract}
Bursting oscillations are commonly seen as a mechanism for information coding in neuroscience and have also been observed in many physical, biochemical, and chemical systems. This study focuses on the computational investigation of mixed-mode bursting oscillations (MMBOs) generated by a simple two-dimensional integrate-and-fire-or-burst (IFB) model. We demonstrate a new paradigm for the generation of MMBOs, where birhythmicity and noise are the key components. In the absence of noise, the proposed model exhibits birhythmicity of two independent bursting patterns, bursts of two spikes and bursts of three spikes, depending on the initial condition of the model. Noise induces the random transitions between two bursting states which leads to MMBOs, and the transition rate increases with the noise intensity. Our results provide a systematic view of the roles of noise and initial condition: the bursting dynamics produced by the proposed model heavily rely on the initial conditions when noise is weak; while for intermediate and strong noise, the burst dynamics are independent of the initial condition.
\end{abstract}

\maketitle

%\begin{quotation}
%\end{quotation}

\section{Introduction}

Bursting oscillations are an essential rhythm of biophysical activity, in particular in the brain regions\cite{Zeldenrust2018}, and have also been observed in many other systems such as electronic circuits \cite{Wijekoon2008,Savino2009}, lasers\cite{Kerse2016,Ruschel2017}, chemical reactions \cite{Tian2013,Beims2018}, and gene expression \cite{Nicolas2018}. The pattern of bursts is characterized by the periodic transition between a quiescent phase (without spikes) and an active phase (with two or more spikes with high firing rates). Thus bursting oscillations have a slow-fast dynamics formalism. They have been classified as various types based on mathematical criteria, particularly bifurcation theory \cite{Izhikevich2000, Desroches2020TowardsAN}.

This paper focuses on a relatively new type of burst, mixed-mode bursting oscillations (MMBOs). The term MMBOs was introduced as an analogy to mixed-mode oscillations (MMOs) \cite{Desroches2013MMBO} and the latter represent an alternation between small-amplitude oscillations (SAOs) and large-amplitude oscillations (LAOs)\cite{Brons2008}. But in MMBOs, repetitive fast spikes, instead of single spikes, occur within LAOs. Several mechanisms have been reported to contribute to the generation of MMBOs. Ref.~\onlinecite{Desroches2013MMBO} showed that MMBOs appear in dynamical systems with a spike-adding mechanism where saddle-type canards allow one or more LAOs in each burst event. Ref.~ \onlinecite{Rubin2015} reported that, for a single oscillator with one fast and one slow variable, simple MMOs become MMBOs when its adaptation changes. Ref.~\onlinecite{Ghosh2020} found a similar change from MMOs to MMBOs for a random network with diverse excitable neurons when the coupling is at intermediate level. Ref.~\onlinecite{Leutcho2020MMBO} reported that MMBOs could also be caused by bubbles of bifurcation. Ref.~\onlinecite{Bacak2016} revealed that heterogeneity (e.g. distributed neuronal excitability, diverse cellular properties, and network connections) is the key to produce MMBOs in a heterogeneous network.  

In this work, we introduce and explain a new mechanism for the generation of MMBOs using a simple two-dimensional bursting model. The key elements in the paradigm are birhythmicity and noise. Birhythmicity refers to the coexistence of two stable limit cycles \cite{Decroly1982} and they can be tonic spiking and bursting\cite{Shilnikov2005}, LAOs and SAOs \cite{Desroches2012}, or even chaotic oscillators\cite{Pisarchik2006}. The first main new result of our study is the birhythmicity of two independent burst patterns (bursting oscillations with two spikes and three spikes) depending on the initial conditions of the model. Birhythmicity is a phenomenon existing in many systems in various disciplines \cite{Wiehl2021}; however, to our knowledge, our work is the first computational study to demonstrate the birhythmicity of two periodic bursting patterns. The second main new result is that noise induces the transitions of two rhythmic activities (i.e. two burst patterns) and consequently result in MMBOs. This mechanism is commonly known as noise-induced transitions \cite{Bashkirtseva2015}. We then divide the noise intensity into four categories (weak, intermediate, strong, and extra-strong) and study how the initial conditions and noise intensity affects bursting dynamics. The results of this article may be generalized to other problems where a birhythmicity of two burst modes exists. 

The paper is organized as follows. In Sec.~II, we introduce the neuronal models. In Sec.~III, we firstly demonstrate the birhythmicity of two independent bursting modes in the absence of noise, then illustrate that noise evokes the switch between two bursting modes, and consequently, MMBOs are generated. We characterize the influence of noise on the dynamics of MMBOs for the trails with two different initial conditions, followed by the systemic view of the dynamics of MMBOs governed by two factors, noise, and initial conditions. An extra-strong noise is used to further explain the impact of noise on MMBOs. The paper ends with a conclusion and discussion in Sec.~IV.

\section{Model}

We adapt the integrate-and-fire-or-burst (IFB) model\cite{Smith2000} for its simplicity and widely acceptance as an important prototype for bursting. Similar to the classic integrate-and-fire (IF) model, IFB has two variables, $v$ for membrane potential of the neuron, and $h$ for the inactivation gating variable of the calcium current. Both IF and IFB models rely on a firing threshold, $v_\theta$, and a reset membrane potential, $v_{reset}$. Thus, when $v_t \geq v_\theta$, $v_{t+\Delta t} = v_{reset}$ where $\Delta t$ is the time step. We introduce an additive noise, $D\xi(t)$, to account for the local noise such as environmental fluctuations. Therefore, the model for our simulations was governed by the following differential equations: 
\begin{align}
C\frac{dv}{dt} &= I_0 + I_1 \cos(2\pi ft) - I_L - I_T + D\xi  \\ 
\frac{dh}{dt} &= \begin{cases}
                 -h / \tau_h^- & \text{if $v<v_h$}\\
                 (1 - h) / \tau_h^+  & \text{if $v>v_h$}
                  \end{cases}
\end{align}
$I_L$ is the leakage current, $I_L$ = $g_L$($v - v_L$), where $g_L$ and $v_L$ are the conductivity and reversal potential for $I_L$, respectively. $I_T$ is the $Ca^{2+}$ current, $I_T$ = $g_T m_{\infty}h(v - v_T)$. The voltage-dependent activation gating, $m_{\infty}$, is described by the Heaviside step function $m_{\infty}$ = $H(v - v_h)$, and $g_T$ is the conductivity for $I_T$. The additive noise has two compartments where $D$ is the noise intensity and $\xi(t)$ is Gaussian white noise with a mean of zero and standard deviation of 1. The parameter values are listed in Table 1.

The gating variables, $m_{\infty}$ and $h$, define the bursting behaviour of a neuron. Whenever $v < v_h$, $v$ is hyperpolarized and $h$ starts to increase towards 1 with the time constant $\tau_h^-$. This causes the calcium current to be deinactivated and therefore, bursting behaviour is observed. On the other hand, when $v > v_h$, $v$ is not hyperpolarized and $h$ starts to decrease to 0 with another time constant $\tau_h^+$, thus causing the calcium current to be inactivated. Due to the unequal values of $\tau_h^-$ and $\tau_h^+$, this system actually has three time scales, a fast time scale for $v$, and a slow and an extra slow time scales for $h$.

The IFB model itself relies heavily on initial conditions to determine the firing mode of a neuron. However, introducing a noise signal, $D\xi$, in the voltage differential equation helps us analyze the behaviour of the membrane potential aside from its deterministic case. By changing the value of $D$, it is possible to look into the effects of the noise in the spikes per burst, as well as how significant the initial conditions for $v$ and $h$ are as $D$ changes. Using $D = 0$ gives us the deterministic case, while using $D > 0$ gives us stochastic cases. 

MATLAB is used to simulate the model and perform numerical analysis. The Euler-Maruyama method is used to approximate the numerical solutions to Equs.~(1)-(2) with the time step $\Delta t$ = 1/50 ms. The results presented in Figs.~2 and 6 are averaged over 40s for one trial, and the results presented in Figs.~4-5 averaged over 30s for 300 trials.

\begin{table}
\caption{Parameter values}
\begin{center}
\begin{tabular}{ccc|ccc}
\hline Parameter & Value & Unit & Parameter & Value & Unit\\\hline
$C$ & 2 & $\mu$F & $v_L$ & -65 & mV \\
$v_h$ & -60 & mV & $v_T$& 120 & mV \\
$v_{\theta}$ & -35 & mV & $v_{reset}$& -50 & mV \\
$g_L$ & 0.035  & mS &$g_T$ & 0.07 & mS\\
$f$ & 5 & Hz & $I_0$ & -0.05 & $\mu$A\\
$I_1$ & 1.6 & $\mu$A &  $\tau_h^+$ & 200 & ms\\
$\tau_h^-$ & 20 & ms &  &  & \\
\hline
\end{tabular}
\end{center}
\end{table}

\section{Results}

\subsection{Birhythmicity of deterministic bursting modes}

The deterministic system of Equs.~(1)-(2) exhibits two independent bursting patterns, depending on the initial conditions of $v$ and $h$, denoted as ($v_0$, $h_0$). To show the sensitivity of this deterministic system to the initial conditions, we use two pairs of initial conditions with a tiny difference, $(v_0, h_0)$=(-45, 0.045) and (-45, 0.05). Bursts with two spikes (blue trace in Fig.~1a) and bursts with three spikes (yellow trace in Fig.~1a) are produced respectively. The blue and yellow bars on the top of the time series of the membrane potential label the spike times and they are color-coded with the time series. For the sake of simplicity, we call the 2-spike burst "mode 2" and the 3-spike burst "mode 3". In the following sections, in the presence of noise, the system could produce isolated single spikes and bursts with more than 3 spikes and they are denoted as "mode 1" and "mode 4", respectively. The $v$-$h$ phase plane in Fig.~1b provides a better view of this "butterfly effect": a small change on the initial conditions (labeled by solid dots in Fig.~1b) leads to a large difference in the bursting patterns. The blue trajectory has two spikes in one burst cycle, while the yellow trajectory has three spikes per burst cycle.    

After a transient time, the system reaches a steady state. The steady states of these two trajectories are plotted in Fig.~1c, after removing their transient phases. We separate the phase plane into three regions to examine the temporal change of both variables. Regions I and II correspond to a physiological slow subsystem while region III corresponds to a fast subsystem. In region I, both trajectories are nearly overlapped. In region II, the trajectories start to differ when their minimal membrane potentials ($v_{min}$) are reached ($v_{min}\approx$ -87 mV for the blue trace and -89 mV for the yellow trace). The blue trajectory corresponding to burst mode 2 is above the yellow trajectory (i.e. burst mode 3). Then $h$ reaches it maximal value ($h_{max}$) when $V$=-60mV in region II (i.e. the end of slow-subsystem). $h_{max}$ corresponding to mode 2 is smaller than $h_{max}$ from mode 3 ($h_{max}\approx$ 0.42 in mode 2 and 0.44 in mode 3). This further causes the different numbers of spike firing in region III (the fast subsystem). The blue trajectory with a smaller $h_{max}$ has a less number of spikes than the yellow trajectory with a larger $h_{max}$ as shown in region III of Fig.~1c.

\begin{figure}
\includegraphics[width=3.37in]{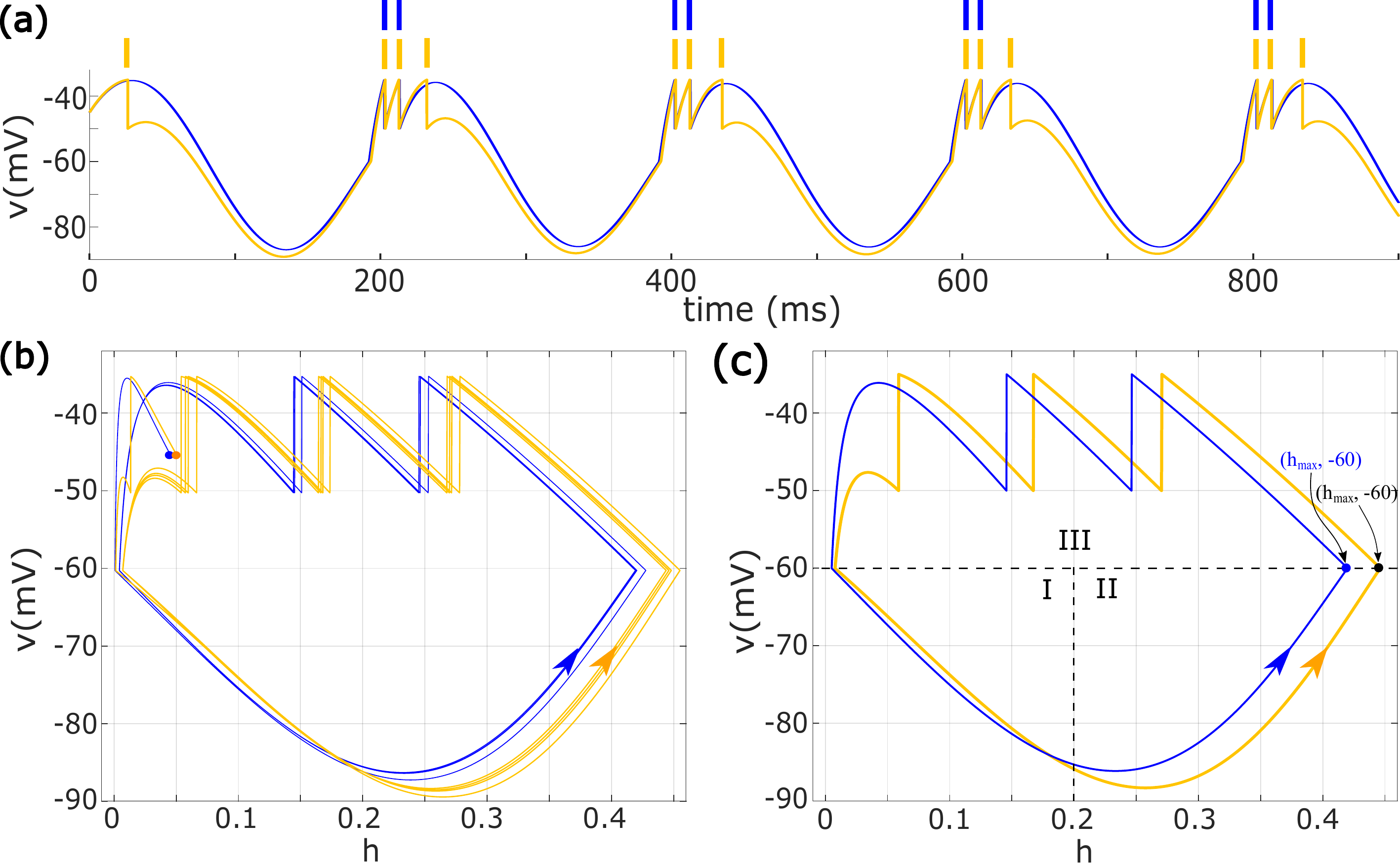}
\caption{\label{fig1} (a) The blue and yellow traces of the membrane potential result from the deterministic system of Equ.~1 with two initial conditions, $(v_0, h_0)$=(-45, 0.045) and (-45, 0.05), respectively. The spike times are marked by the vertical bars on the top of the voltage traces. The voltage traces and spike times are color-coded: 2-spike bursts (burst mode 2) in blue and 3-spike bursts (burst mode 3) in yellow. (b) The corresponding trajectories were plotted on the phase plane over the time course from 0 to 900ms. The colored dots label the initial conditions and the arrows indicate the direction of the trajectories. (c) The steady-state trajectories, after the transient phase is removed. The slow subsystem of the model is in regions I and II, and the fast subsystem is in region III. The blue and black solid circles label the location when $h$ reaches its maximum value, $h_{max}$, where $v$=-60 mV.}
\end{figure}

To get a thorough understanding of birhythmicity of this system, we compute the average number of spikes per burst with all possible pairs of $v_0\in$[-90, -35] and $h_0\in$[0,1] as the initial conditions. The results are presented in Fig.~2, and it illustrates that the deterministic system has two stable bursting modes: mode 2 (blue) and mode 3 (yellow). For the majority of the initial condition pairs, the system results in burst mode 2. Surprisingly, the initial values pairs leading to mode 3 form multiple tilted strips (narrow yellow bands in Fig.~2). Moreover, both the width and length of these tilted strips get smaller as $h_0$ moves from 0 to 1. For two pairs of initial conditions used in Fig.~1, $(v_0, h_0)$=(-45, 0.045) is on the blue region and $(v_0, h_0)$=(-45, 0.05) is in the yellow region.

\begin{figure}
\includegraphics[width=3in]{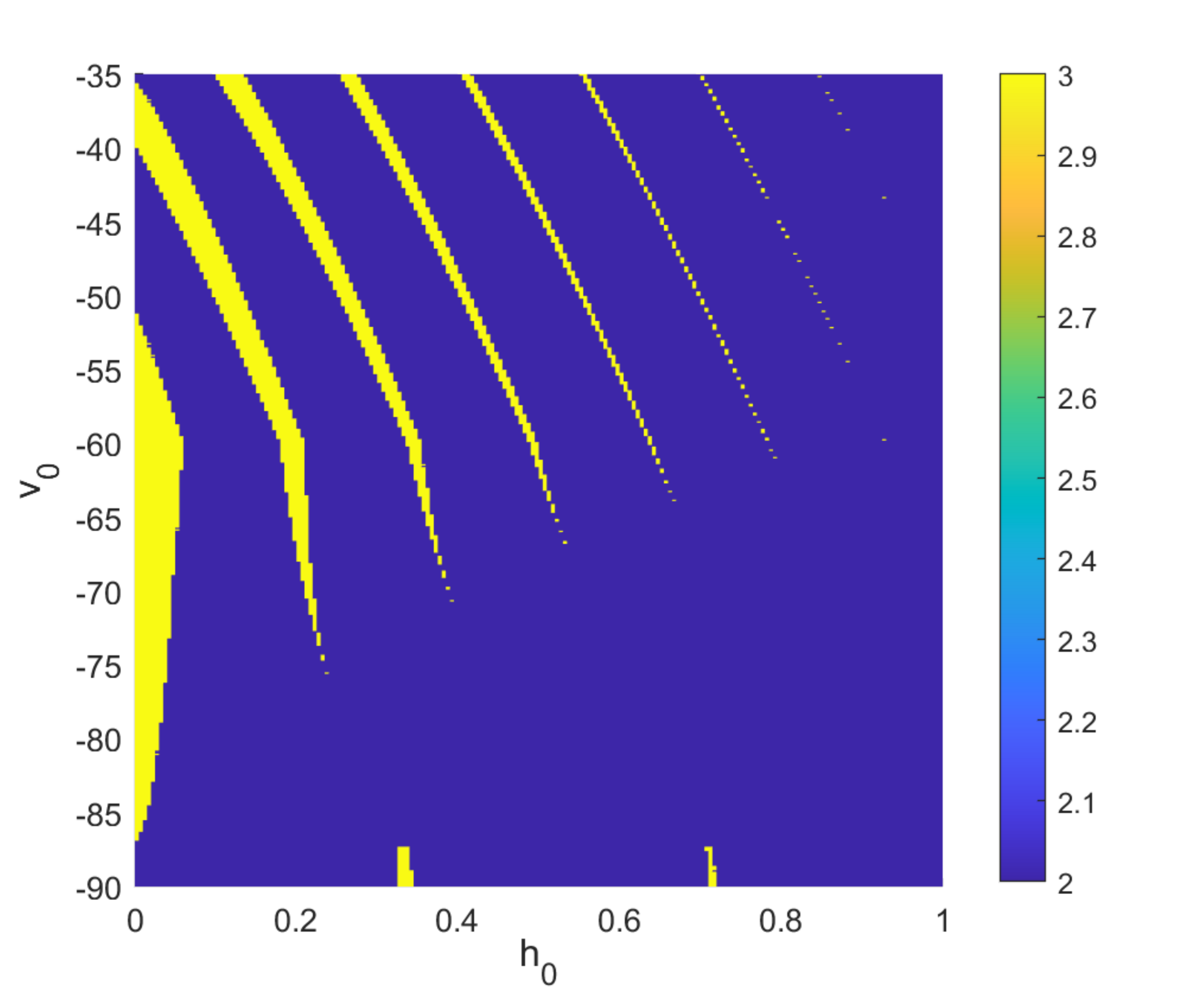}
\caption{\label{fig2} Partitions of ($v_0$, $h_0$) space (for fixed values of the other parameters) according to the number of spikes per burst (color-coded by the color bar on the right) of the deterministic system ($D=0$). This diagram shows two bursting modes: burst mode 2 (blue) and burst mode 3 (yellow) for $v_0\in[-90,-35]$ and $h_0\in[0,1]$. }
\end{figure}

\subsection{MMBOs induced by noise and birhythmicity}

Noise is then introduced in the system. It triggers a random transition among different bursting modes over time, thus noise-induced MMBOs are formed as illustrated in Fig.~3. Different noise intensities impact the bursting dynamics in different ways. Three example $D$ values (0.06, 0.5, and 2) in Fig.~3 represent weak, intermediate, and strong noise intensities respectively. Three observations are made from Fig.~3. Firstly, for weak and intermediate noise (top two rows of Fig.~3a), there are only two burst modes (mode 2 and mode 3). However, for strong noise (bottom row of Fig.~3a), more burst modes appear (mode 1 and mode 4). Secondly, the switching between different burst modes is more frequent as $D$ increases. As shown in Fig.~3a, for $t\in$ [100ms, 3000ms] (i.e. the first isolated spike is truncated), the burst mode switches only once with $D$ = 0.06; there are 5 switches between two burst modes (modes 2 and 3) with $D$ = 0.5, and 9 switches between four burst modes (modes 1 to 4) $D$ = 2. Thirdly, noise has more influence on the slow subsystem (i.e. regions I and II where $v$ < 60mV), thus $h_{max}$ of each burst cycle when $v$ = -60mV in region II (i.e. the end of the slow subsystem) is more diversely distributed with a higher $D$, which further causes the genesis of different burst modes (see the number of firing in region III). Similar to the deterministic case demonstrated in Fig.~1c, a lower $h_{max}$ leads to a burst mode with a fewer number of spikes. 

\begin{figure}
\includegraphics[width=3.37in]{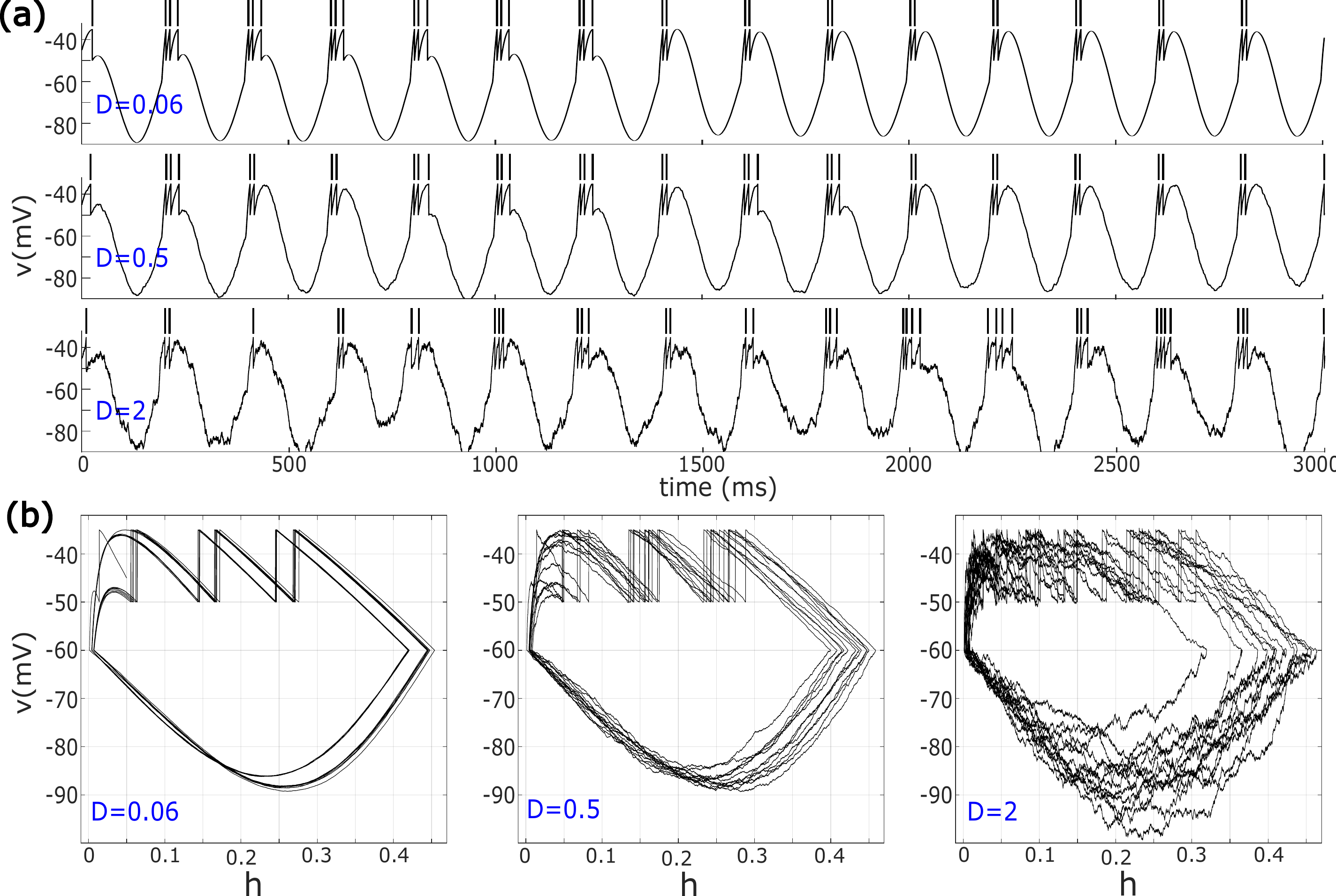}
\caption{All panels have the same initial condition $(v_0, h_0)$=(-45, 0.05), corresponds to deterministic burst mode 2 (i.e. 3 spikes/burst) in Fig.~1. (a) The example voltage traces with D=0.06 (top, weak noise), 0.5 (middle, intermediate noise), and 2 (bottom, strong noise). The spike times are marked by the vertical bars on the top of the voltage traces. (b) The corresponding trajectories in the phase planes.}
\end{figure}

\subsection{The roles of noise on bursting dynamics}

We compute the transition rates (the number of switches between different burst modes per second) and the occurrence percentages of each burst mode (the number of one burst mode over a total number of all bursts in one trial) averaged over 300 trials of 30-second/trial starting from each initial condition, $(v_0, h_0)$=(-45, 0.045) and (-45, 0.05). They are the functions of noise intensity and are plotted in Fig.~4. Based on these statistical results, we categorize noise into four categories: weak (0<$D\leq$0.14), intermediate (0.14<$D\leq$1.2), strong (1.2<$D\leq$5), and extra-strong ($D$>5). The first threshold $D$=0.14 is chosen because the transition rates curves (Fig.~4ab) and the occurrence percentage curves corresponding to one burst mode (Fig.~4bc) from both initial conditions are identical for $D\geq$0.14. The second threshold $D$=1.2 is chosen because burst modes 1 and 4 starts to occur when $D\geq$1.2. The reason to choose the third threshold $D$=5 will be explained in section III.F.

The transition rate in Fig.~4ab increases with a larger $D$, which agrees with our observation on the example trials in Fig.~3. The transition rate corresponding to $h_0$=0.045 is very close to zero over the weak noise regime, which implies that the burst mode remains almost the same as the deterministic case (i.e. barely switch to other modes) for most of the weak noise trials. But for $h_0$=0.05, the transition rate between mode 3 (deterministic case) and mode 2 is still very low over the weak noise regime, only a few times over each 30-second trial (see Fig.~3a top for an example over 3-second). Therefore the transition rate corresponding to $h_0$=0.05 is slightly higher than that of $h_0$=0.045 over $0.03\leq D\leq 0.14$. For intermediate and strong noise (more exactly, $D>0.14$), the average transition rates result from two initial values are almost identical.

\begin{figure*}
\includegraphics[width=5.5in]{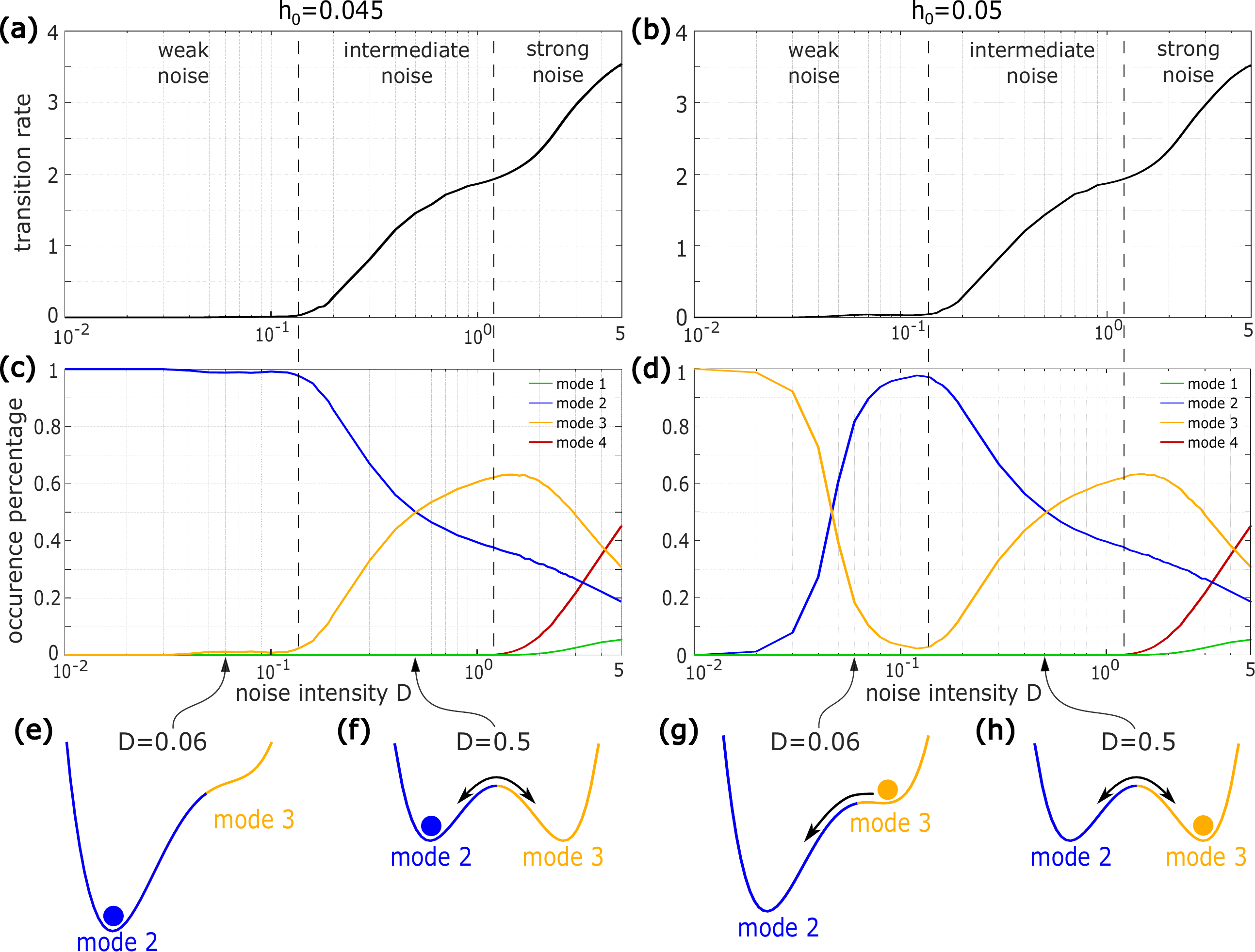}
\caption{Side-by-side comparison of burst dynamics with two initial conditions, $(v_0, h_0)$=(-45, 0.045) on the left and (-45, 0.05) on the right and weak (0<$D$<0.14), intermediate (0.14<$D$<1.2) and strong noise (1.2<$D$<5). (a)(b) The transition rate (the number of switches among different burst modes per second) v.s. noise intensity $D$. Two transition rate curves are almost identical. (c)(d) The occurrence percentages of four burst modes v.s. noise intensity $D$. For one $D$ value, the sum of occurrence percentages of all four modes is 1. (e)-(h) The "double-well" diagrams to illustrate the switches between modes 2 and 3 for $D$=0.06 (weak noise intensity) and $D$=0.5 (intermediate noise intensity). The blue balls are located at mode 2 because it is the deterministic mode determined by the initial condition, $(v_0, h_0)$=(-45, 0.045). The yellow balls are in mode 3 which is the corresponding deterministic burst mode for the initial condition, $(v_0, h_0)$=(-45, 0.05). All curves in (a)-(d) are averaged over 300 trials, 30 seconds per trial.}
\end{figure*}

In subsection III.A, we have shown that the deterministic burst modes depend on the initial conditions. Initial conditions still produce a major impact on the burst dynamics, especially over the weak noise range (0<$D$<0.14). With $(v_0, h_0)$=(-45, 0.045), corresponding to deterministic burst mode 2, the average occurrence percentages of mode 2 (blue curve in Fig.~4c) decreases very slowly from 100\% to 93\% over the weak noise regime, then decreases much faster from 93\% to 19\% over the intermediate and strong noise regimes. Meanwhile, the average occurrence percentage of mode 3 (yellow curve in Fig.~4c) increases very slowly over the weak noise regime then rapidly over the intermediate noise regime. It reaches the highest value, 63\%, at $D$=1.5 and then decreases because the strong noise induces another two burst modes (modes 1 and 4). With $(v_0, h_0)$=(-45, 0.05), corresponding to deterministic burst mode 3, the average occurrence percentages of mode 3 (yellow curve in Fig.~4d) dramatically decreases from 100\% to 2\% for 0<D<0.12 (weak noise), then increases over the intermediate noise regime. Correspondingly, the percentage curve of mode 2 (blue curve in Fig.~4d) increases very quickly from 0\% to 98\% for 0<D<0.12. We also notice that, for each burst mode in the stochastic system, its occurrence percentages resulting from two initial conditions are equivalent when $D$>0.14. It means that the initial condition has very limited influence on MMBOs induced by intermediate and strong noise.

The MMBOs induced by weak and intermediate noise due to the birhythmicity of burst modes can be illustrated by the "double-well", which has been widely used to demonstrate the transition between two states of a dynamical system \cite{Aubry1974,Theocharis2006}. As shown in Fig.~4(e-f), each well represents one burst mode, the depth of each well is proportional to the average occurrence percentage of this burst mode and further determines the jump probability from one well to another. The ball indicates the deterministic burst mode determined by the initial condition. Firstly, we consider the role of weak noise. Due to weak noise, the initial burst mode of each stochastic trial is actually the burst mode of this trial in the deterministic case, for example, the initial burst mode of the stochastic trials staring with $(v_0, h_0)$=(-45, 0.05) is mode 3 (see Fig.~3a). For the trails starting with $(v_0, h_0)$=(-45, 0.045), the ball is located at mode 2 initially, and it has a very low probability to jump to mode 3 due to the large difference in the depths of the two wells (Fig.~4e), which leads to a high occurrence percentage of mode 2 (around 98\% for $D$=0.06 in Fig.~4c). For the trails starting with $(v_0, h_0)$=(-45, 0.05), the ball is initially located at mode 3, and it has a high probability to jump to mode 2. Once it jumps to mode 2, the chances of going back to mode 3 are very low, which causes the dramatic drop of the occurrence percentage curve of mode 2 over the weak noise regime in Fig.~4d. When noise intensity is increased to its intermediate regime, the transition of the ball between two wells is largely determined by the noise and the initial position of the ball has a very low influence. For example, when $D$=0.5, the ball has nearly equal chances to jump between two wells for both initial conditions (Fig.~4fh). This explains why the occurrence percentage curves are identical for each burst mode over the intermediate and strong noise regime in Fig.~4cd.

\subsection{Inter-spike intervals of noise-induced MMBOs}

In the absence of noise, the inter-spike interval histogram (ISIH) associated with $(v_0, h_0)$=(-45, 0.045) has two peaks (blue peaks in Fig.~5a) because this initial condition leads to bursts of two spikes (see Fig.~1a). So these two ISIH peaks have the same height and they are located at 11 ms (the interval between two spikes within a burst) and 189 ms (the interval between the last spike of a burst and the first spike of the next burst). The ISIH associated with $(v_0, h_0)$=(-45, 0.05) has three equal-height peaks at 10 ms, 21 ms, and 169 ms (yellow peaks in Fig.~5a). They are the inter-spike intervals (ISIs) between any two successive spikes in burst mode 3.

When noise is weak, the locations of two blue ISIH peaks associated with $(v_0, h_0)$=(-45, 0.045) do not change (see blue peaks in Fig.~5b,c,d), although two additional peaks with very low height are presented at 21 ms and 169 ms when $D=0.1$. It means that the majority of bursts associated with $(v_0, h_0)$=(-45, 0.045) have two spikes per burst (burst mode 2) over the weak noise regime, which agrees with the results in Fig.~4c. But the height and width of these two ISIH peaks decreases and increases, respectively, when $D$ changes from 0.03 to 0.1, because the introduction of weak noise perturbs the timing of deterministic spikes. 

Weak noise leads to more changes on the ISIH associated with $(v_0, h_0)$=(-45, 0.05). When $D$ changes from 0 to 0.03, the number of yellow peaks changes from three to four. The height of the first yellow peak at 10 ms does not change, but the heights on the second and third peaks are lower because the third spike of a small portion of bursts is eliminated by noise (i.e. a small portion of bursts change from mode 3 to mode 2), which agrees with the results in Fig.~4d. The eliminated 3rd spike also causes the generation of a new ISIH peak (i.e. the fourth yellow peak at 189 ms), although this peak has a very low height as shown in Fig.~5b. When $D$ keeps increasing (e.g. $D$=0.06 and 0.1), the heights of second and third ISIH peaks are tremendously decreased but the heights of the first and fourth peaks are relatively higher, which implies that the majority of bursts change to mode 2, as we have seen in Fig.~4d. 

When $D=0.1$, the blue peaks and yellow peaks are almost identical, and when $D>0.14$ the ISIHs associated with both initial conditions are identical. It implies that the initial condition has less and less influence on burst dynamics when $D$ increases from weak to intermediate or strong regime, same conclusion we have drawn in section III.C. We also notice that, as $D$ increases from 0.3 to 1 (over the strong noise regime), the first two peaks gradually merge into one peak with a much larger width (Fig.~5ef), and similar behavior for the third and fourth peaks, which indicates that strong noise is able to evoke all four burst modes.

\begin{figure*}
\includegraphics[width=6in]{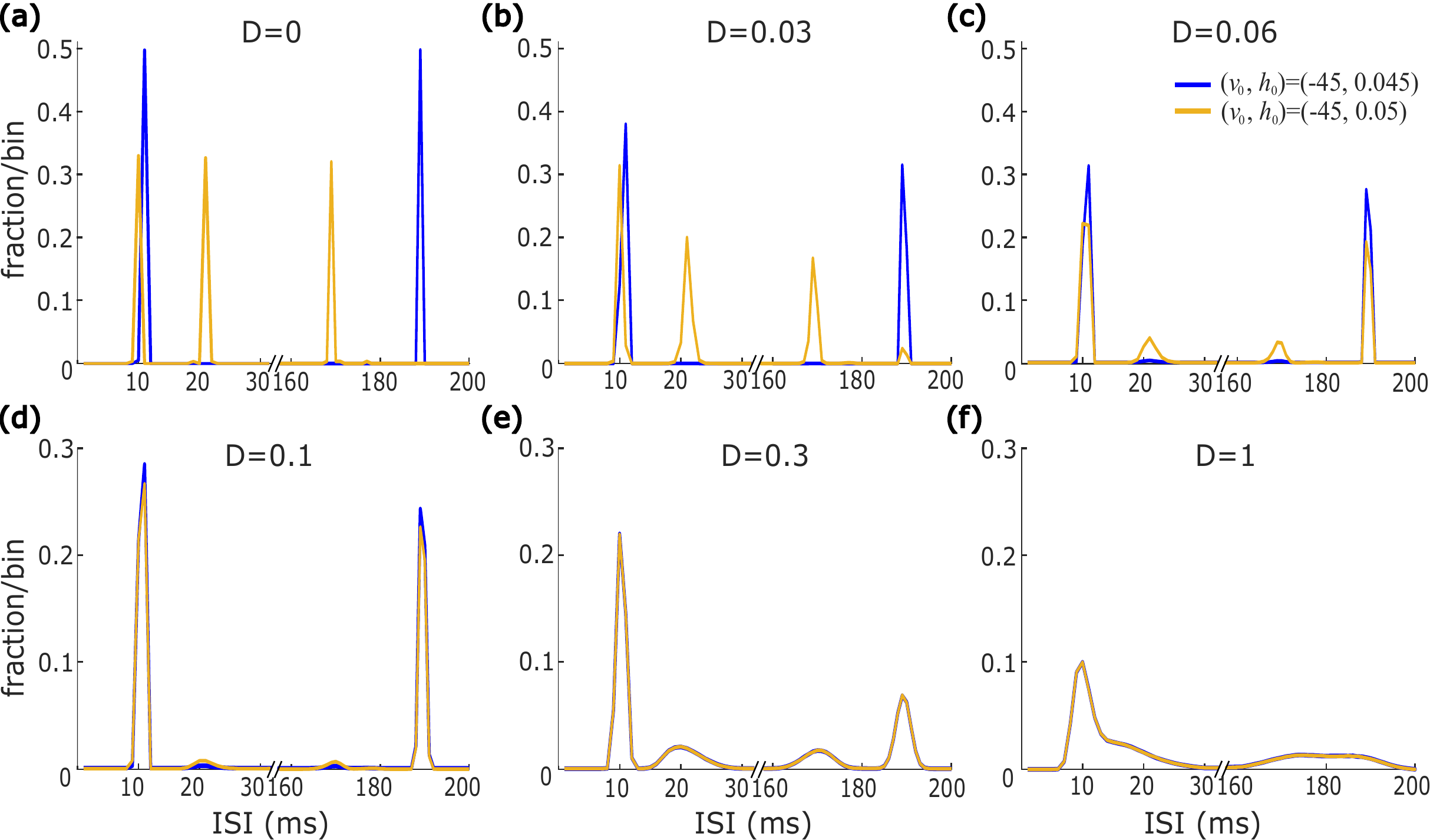}
\caption{Inter-spike interval histograms (ISIHs): fraction/bin v.s. inter-spike interval (ISI) for various noise intensities, $D=$ (a) 0, (b) 0.03, (c) 0.06, (d) 0.1, (e) 0.3, and (f) 1. The binwidth of ISI is 1ms. Blue and yellow traces are associated with $(v_0, h_0)$=(-45, 0.045) and (-45, 0.05), respectively. The results are averaged over 300 trials (30-second/trial).}
\end{figure*}

\subsection{The roles of noise and initial condition}

To characterize the transition from birhythmicity to MMBOs regulated by both initial condition and noise, we compute the average number of spikes per burst for all possible ($v_0$, $h_0$) $\in$ [-90, -35]$\times$[0,1] and four noise intensities ($D$=0.025, 0.05, 0.1, and 0.5). The results averaged over 40 seconds are presented in Fig.~6, one trial for each initial condition pair.

Compared with Fig.~2 ($D$=0), the diagram with very weak noise ($D$=0.025) in Fig.~6a  does not have obvious changes except that the edges of the narrow yellow strips are vague. It is in line with the results in Fig.~4, where the average occurrence percentage of burst modes 2 is 100\% for $(v_0, h_0)$=(-45, 0.045) but the average occurrence percentage of mode 3 decreases to around 95\% for $(v_0, h_0)$=(-45, 0.05) and the latter initial condition is on the edge of yellow stripes. When $D$ increases to 0.05 (Fig.~6b), the blue region associated with burst mode 2 barely changes (the average occurrence percentage of mode 2 is around 99\%). However, the narrow yellow stripes associated with burst mode 3 become unstable due to trial-to-trial variation of MMBOs (the averaged occurrence percentage of mode 3 is around 39\%). When $D$ is close to the threshold between weak and intermediate noise ($D=0.1$, Fig.~6c), almost all bursts are mode 2 regardless of the initial condition, which matches with the average occurrence percentages of mode 2 (99\% in Fig.~4c and 96\% in Fig.~4d). For intermediate or strong noise, the average number of spikes per burst is uniform on the plane of $v_0$ vs. $h_0$ (Fig.~6d). 

These results indicate that the initial condition has an important influence on the burst patterns if noise is zero or relatively small (Fig.~6ab), and the roles of initial condition and noise cause various burst dynamics. However, the initial condition has very small or even no impact on MMBOs induced by relatively larger noise (Fig.~6cd).

\begin{figure*}
\includegraphics[width=6in]{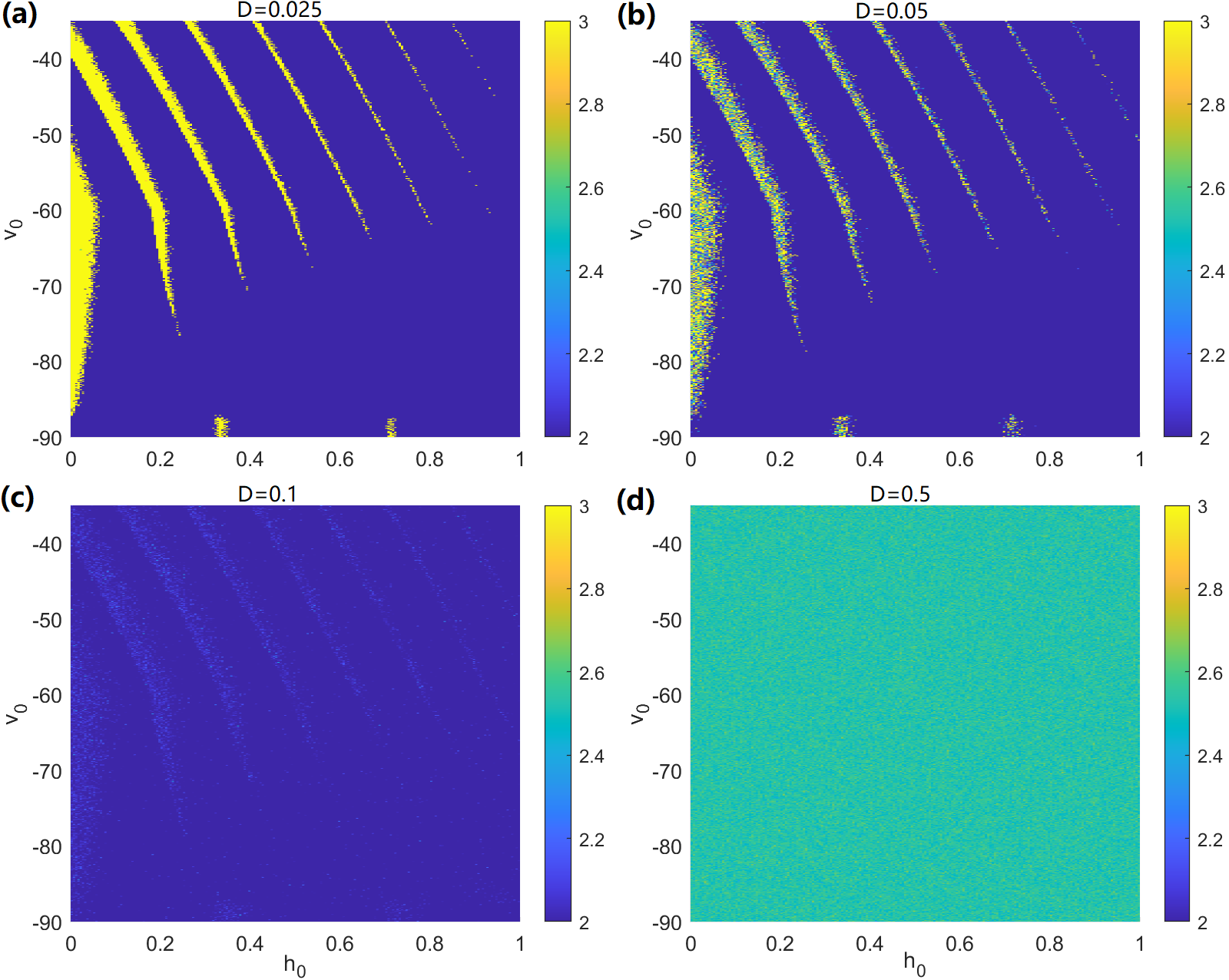}
\caption{The average number of spikes per burst (color-coded by the color bar on the right) of the stochastic system starting from each pair of initial conditions, ($v_0$, $h_0$) $\in$ [-90,-35]$\times$[0,1], but with different noise intensities: (a) $D$=0.025, (b) $D$=0.05, (c) $D$=0.1, (d) $D$=0.5. The results presented here are averaged over 40 seconds, one trial for each initial condition pair.}
\end{figure*}

\subsection{MMBOs driven by extra-strong noise}

Fig.~7a demonstrates the ISIH curves for one strong noise $D=3$ and three extra-strong noise, $D=$ 5, 7, and 10. Unlike most ISIH curves in Fig.~6, ISIH curves here have two peaks and each peak has a much larger peak width. The first peak at the shorter ISI corresponds to the intra-burst interspike intervals, and the second peak at longer ISI corresponds to the intervals between the last spike of a burst and the first spike of the next burst. $D=5$ is chosen as the threshold between strong and extra strong noise regimes because strong noise ($1.2<D<5$) has two isolated ISIH peaks while extra-strong noise ($D>5$) has two connected peaks, that is, the "trough" (i.e. the minimum fraction/bin) in the middle of two ISIH peaks is greater than zero as shown in Fig.~7a. The ISI value at the "trough" is often used to identify the occurrence of a burst \cite{Selinger2007}, and here 80 ms is used to differentiate burst spikes from intra-burst spikes. But extra-strong noise dominates the system (Fig.~7c) and drives the system to randomly produce spikes. Thus ISIs are widely distributed when $D>5$ extra-strong, instead of concentrated on two peaks. As a result, it is difficult to determine if some spike, for example, the spike labeled by the blue vertical bar at around 1900 ms in Fig.~7b, belongs to the burst prior to it, or the burst after it, or if it is an independent single spike (i.e. mode 1). We also notice that the second ISIH peak starts to disappear when $D\geq10$, which implies that, with continuously increased noise over the extra-strong regime, the system is forced to generate more random individual spikes, instead of bursts of spikes.

\begin{figure}
\includegraphics[width=3.37in]{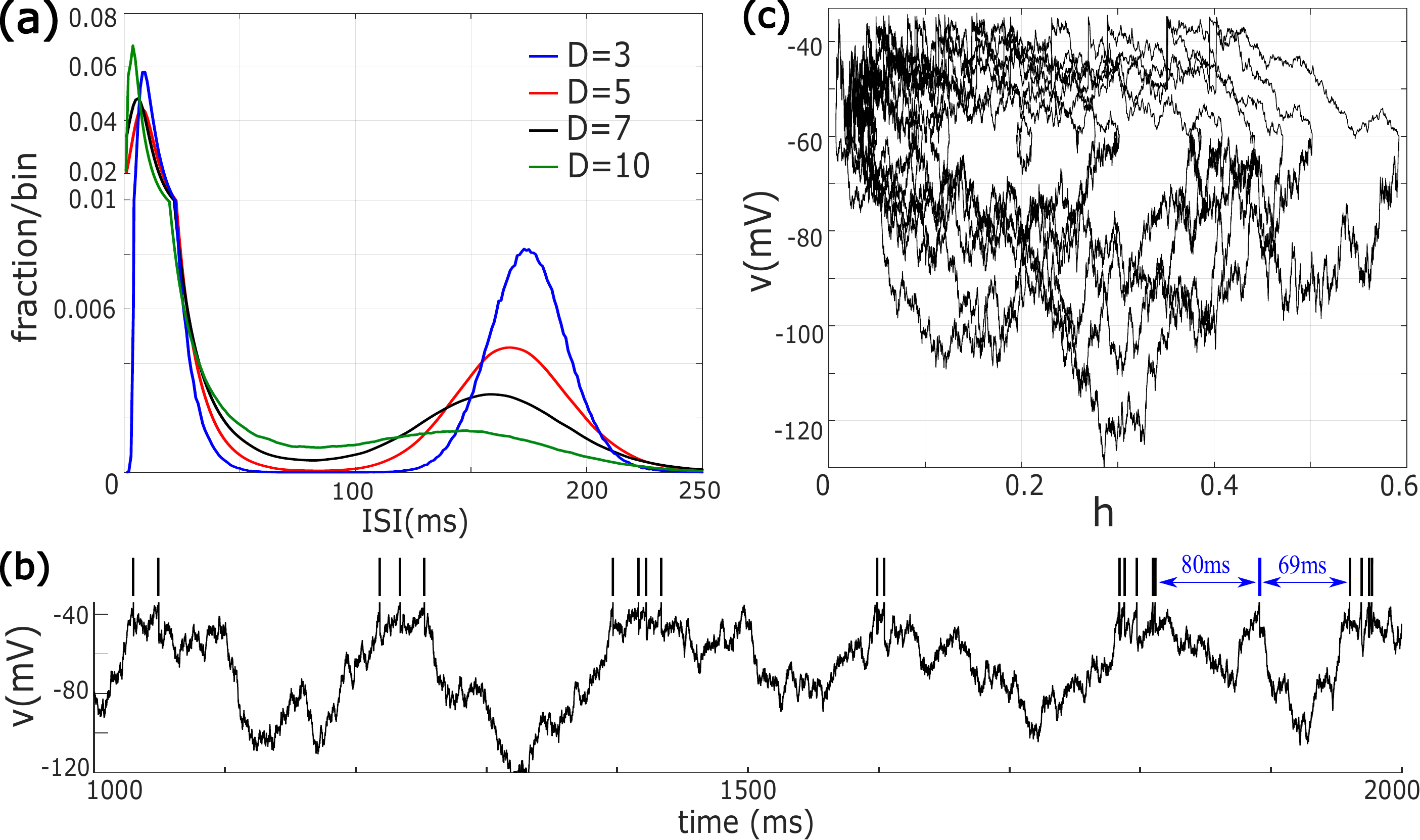}
\caption{(a) The ISIH curves for four noise intensities $D=$3, 5, 7, 10. The binwidth of ISI is 1 ms. Four ISIH curves are identical for ($v_0$, $h_0$)=(-45, 0.045) and (-45, 0.05). (b) An example time series of membrane potential with $D=7$ and ($v_0$, $h_0$)=(-45, 0.045). The spike times are marked by the vertical bars on the top of the voltage traces. (c) The corresponding trajectories were plotted on the phase plane. }
\end{figure}

\section{Conclusion and Discussion}

This study presented a new mechanism for MMBOs, where birhythmicity and noise are key factors and performed a systematic examination of how noise and initial conditions affect the dynamics of MMBOs. We concluded our major findings as to the following: 
\begin{itemize}
\item In the absence of noise, the model shows the birhythmicity of two independent burst patterns depending on the initial conditions (Figs.~1-2).
\item when noise is weak, the initial condition has a major influence on the burst patterns. For initial conditions associated with deterministic burst mode 2, weak noise barely modifies the burst modes (blue region in Fig.~6abc), while for initial conditions associated with deterministic mode 3, weak noise triggers the interaction between two burst modes and the dynamics of burst patterns varies trial by trial (color-mixed stripes in Fig.~6b).
\item when noise intensity is close to the threshold between weak and intermediate noise ($D\approx$ 0.14), almost all bursts are mode 2 whatever the initial condition is (Fig.~5c).
\item When noise intensity is intermediate or strong, the time evolution of MMBOs is independent of initial values but greatly influenced by noise (Fig.~4cd, Fig.~5ef, and Fig.~6d) because of the increased transition rate between different burst modes (Fig.~4ab).
\item Strong and extra-strong noise also triggers other burst modes (i.e. bursts with one or more than three spikes).
\end{itemize}

Although the focus of this work is the noise-induced MMBOs, we presented a preliminary computational result of the birhythmicity of two independent burst patterns using a simple 2-dimensional deterministic IFB model. Thus, with a dynamical systems viewpoint, it will be essential to analyze the nonlinear properties of this deterministic system (e.g. identifying the parameter regions for the LFB model where the birhythmicity of two burst modes exists) in the near future. Additionally, there are a number of challenges to consider, such as investigating if such birhythmicity occurs in other bursting models for example adaptive exponential integrate-and-fire (AdEx) model \cite{Brette2005} and Hindmarsh-Rose (HR) model\cite{Hindmarsh1984}, and generalizing the results from LFB model to other problems where such a birhythmicity exists.

This study showed that, regardless of the initial conditions, the simple IFB model has the identical characteristics of MMBOs when subject to intermediate or strong noise, which agrees with many numerical and experimental studies that noise facilitates the transition between coexisting solutions \cite{Paydarfar2006, Slepukhina2020}. However, for weak noise, the burst dynamics have a large difference depending on the choice of the initial condition. It implies that one should be very cautious about the influence of initial conditions on a system exhibiting MMBOs or MMOs when noise is weak or extremely weak. 

%\begin{acknowledgments}
%We wish to acknowledge the support of the author community in using REV\TeX{}, offering suggestions and encouragement, testing new versions, \dots.
%\end{acknowledgments}

%\section*{Data Availability Statement}

%AIP Publishing believes that all datasets underlying the conclusions of the paper should be available to readers. Authors are encouraged to deposit their datasets in publicly available repositories or present them in the main manuscript. All research articles must include a data availability statement stating where the data can be found. In this section, authors should add the respective statement from the chart below based on the availability of data in their paper.

%\appendix
%\section{Appendixes}

\nocite{*}
%\bibliography{aipsamp}% Produces the bibliography via BibTeX.
\bibliography{ref}

\end{document}